\documentclass[a4paper,fleqn]{cas-sc}

\usepackage[numbers,sort&compress]{natbib}
\usepackage{algorithm}
\usepackage{algorithmicx}
\usepackage{algpseudocode}

\newcommand{\cB}{\mathcal{B}}

\newcommand{\cL}{\mathcal{L}}

\newcommand{\bx}{\mathbf{x}}

\newcommand{\bzero}{\mathbf{0}}

\newcommand{\ti}[1]{\text{\scriptsize{#1}}}
\newtheorem{example}{Example}[section]

\begin{document}
\let\WriteBookmarks\relax
\def\floatpagepagefraction{1}
\def\textpagefraction{.001}

\shorttitle{MPU-PINNs}    

\shortauthors{Yang,Shin,Choi,and Lee}  

\title [mode = title]{A new strategy for physics-informed neural networks based on hierarchical collocation point refinement}      
\tnotemark[1] 
\tnotetext[1]{This work is sponsored by ...}

\author[1]{Minjae Choi}
%\ead{minjaechoi@yonsei.ac.kr}
\affiliation[1]{
    organization={School of Mathematics and Computing, Yonsei University},
    addressline={50 Yonsei-ro, Seodaemun-gu},
    city={Seoul},
    postcode={03722},
    country={Republic of Korea}
}

\author[1]{Dukhwan Shin}
%\ead{2025311235@yonsei.ac.kr}

\author[1]{Youngsoo Yang}
%\ead{***}

\author[1]{Eunjung Lee}[orcid=0000-0001-9989-3555]
\cormark[1]
%\ead{eunjunglee@yonsei.ac.kr}
\cortext[1]{Corresponding author}



\begin{abstract}
Physics-informed neural networks (PINNs) offer a flexible framework for solving partial differential equations (PDEs), but training can become computationally expensive when a large number of collocation points are required to accurately enforce the governing equations. To alleviate this cost, we introduce multigrid-based parameter-updated PINNs (MPU-PINNs), a coarse-to-fine training strategy that progressively increases the number of training points throughout the learning process. 
The proposed approach begins by training a neural network on a coarse set of collocation points and then transfers the learned parameters to successively finer levels. This initialization strategy enables the network to capture the global features of the solution at a relatively low computational cost before refining local details with additional training points. To further improve performance for high-frequency problems, we incorporate a scaling technique that mitigates the effects of spectral bias during training. We evaluate MPU-PINNs on several benchmark PDEs, including two- and three-dimensional Poisson equations, a convection-diffusion-reaction equation, and the Helmholtz equation. Numerical experiments indicate that MPU-PINNs greatly reduce training time while achieving accuracy comparable to that of conventional PINNs and other representative variants such as SA-PINNs and XPINNs. The results further suggest that the proposed coarse-to-fine learning strategy substantially decreases the optimization effort required at finer levels. Overall, MPU-PINNs provide an efficient single-network training framework that enhances the computational efficiency and scalability of PINNs for a broad range of PDE problems.
\end{abstract}

%% Research highlights
%\item 
%\item 
%\item 

% Keywords

\maketitle

%%%%%%%%%%%%%%%%%%% Main text

\section{Introduction}\label{sec:intro}

Partial differential equations (PDEs) play a fundamental role in modeling a wide range of physical phenomena, including fluid flow, heat transfer, wave propagation, and diffusion processes. Consequently, efficient and accurate solution techniques for PDEs are essential in many areas of science and engineering. Classical numerical methods, such as finite difference, finite element, and finite volume methods, have been widely used for this purpose. In recent years, machine-learning-based approaches have emerged as an alternative paradigm for PDE computation. Among them, physics-informed neural networks (PINNs) incorporate the governing equations and boundary conditions directly into the training objective, enabling the approximation of PDE solutions without requiring labeled data \cite{Raissi}. This feature makes PINNs particularly attractive in situations where simulation data are limited or expensive to obtain, and has led to their successful application to both forward and inverse problems.

A central component of PINNs is the set of collocation points used to enforce the PDE residual. In general, increasing the number of collocation points improves the representation of the governing equations and can enhance the accuracy of the learned solution. However, a larger training set also increases the computational cost of evaluating residuals and computing gradients during optimization. As a result, the training cost of PINNs can become substantial, especially for multidimensional problems or problems requiring a large number of collocation points to resolve complex solution features. The importance of collocation points has motivated extensive research on sampling strategies for PINNs. Residual-based adaptive refinement (RAR) and related methods iteratively update the distribution of collocation points according to the current residual error that allocates additional points to regions where the solution is less accurate \cite{RAR,RAD,FI_PINNs}. Other approaches employ adaptive sampling procedures or generative models to construct informative training sets for high-dimensional PDEs \cite{DasPINNs}. These studies indicate that the selection and distribution of collocation points strongly influence both the accuracy and computational efficiency of PINN training.

At the same time, considerable effort has been devoted to improving the optimization and convergence properties of PINNs. Examples include domain-decomposition frameworks such as XPINNs and cPINNs \cite{XPINNs,cPINNs}, which partition the computational domain into smaller subdomains, as well as techniques designed to mitigate spectral bias, balance competing loss terms, and improve stability for challenging PDEs \cite{FF,NTK,Instable_PINN}. While these developments have significantly enhanced the applicability of PINNs, most approaches still train the network on a prescribed set of collocation points from the outset. Consequently, the computational burden associated with large training sets remains an important consideration in practical applications.

While considerable effort has been devoted to improving collocation-point selection and optimization strategies, comparatively less attention has been paid to how information learned from a coarse collocation set can be exploited to accelerate training on a finer one. Since the network parameters obtained from a coarse discretization already contain useful information about the global structure of the solution, they can potentially provide an effective initialization for subsequent training stages. Motivated by this observation, we propose multigrid-based parameter-updated PINNs (MPU-PINNs).
The proposed method first trains a single neural network using a relatively small number of collocation points at a coarse level. The learned parameters are then transferred to the next finer level as the initialization for further training with additional collocation points. By repeating this procedure, the network progressively refines the solution while retaining information acquired at previous levels. In this way, MPU-PINNs combine parameter transfer and progressive collocation-point refinement within a unified framework. The coarse-level training captures the global structure of the solution at a relatively low computational cost, while the finer levels improve local accuracy through increasingly dense collocation sets. Unlike domain-decomposition approaches, the proposed method employs a single neural network throughout the entire training process and does not require multiple local models or interface conditions.

The remainder of this paper is organized as follows. Section~\ref{sec:model} briefly reviews the standard PINNs formulation. Section~\ref{sec:rpoint} introduces the proposed MPU-PINNs framework and presents numerical results for several benchmark PDEs. Section~\ref{sec:high} investigates high-frequency problems and introduces an additional scaling strategy that further improves the performance of MPU-PINNs in the presence of spectral bias. Finally, Section~\ref{sec:conclude} concludes the paper.

\section{Model problem}\label{sec:model}

We consider a set of partial differential equations (PDEs) in a bounded domain $\Omega \subset \mathbb{R}^d$, $d = 2,3$, given as 
\begin{equation}\label{eq:model}
	\mathcal{A}(U) = F \quad \text{in} \quad \Omega,
\end{equation}
subject to the boundary condition $\cB(U)=G$ on $\Gamma:=\partial\Omega$. Here, $\mathcal{A}$ is a differential operator related to the PDEs, $\mathcal{B}$ specifies the boundary conditions, $U$ is the unknown function, and $F$ and $G$ are prescribed functions. The standard PINNs find an approximation of the solution $U$ by using following loss objective function: 
\begin{equation}\label{eq:total_loss}
\frac{\lambda_r}{N_r}\sum_{j=1}^{N_r}\left|\mathcal{A}(U_\theta(\bx_j))-F(\bx_j)\right|^2 \,\, + \,\, \frac{\lambda_b}{N_b}\sum_{j=1}^{N_b}\left|\mathcal{B}(U_\theta(\bx_j))-G(\bx_j)\right|^2,
\end{equation}
where $U_\theta$ denotes the neural network approximation obtained by the trainable parameters $\theta$. The non-negative $\lambda_r$ and $\lambda_b$ weight the PDE residual and boundary loss, respectively. $N_r$ and $N_b$ denote the numbers of points to compute each losses. This basic framework of PINNs has been extensively studied and applied to a wide range of problems over the years. However, it is well-known that the performance of PINNs can vary depending on the problem, and their use is often accompanied by challenges and potential failures (\cite{FI_PINNs} and references therein). One of the main challenges is slow convergence, as optimization methods may struggle to navigate the complex loss landscape efficiently. So, in the following section, we propose a new algorithm designed to accelerate the convergence of PINNs.

\section{Multigrid-based parameter update strategy}\label{sec:rpoint}

We propose a new training strategy that gradually increases the complexity of the optimization problem by progressively enlarging the set of collocation points. When only a small number of collocation points is used, each training epoch is computationally inexpensive because the loss function is evaluated on fewer samples. Consequently, the coarse-level problem can be solved at a lower cost which allows the neural network to efficiently capture the global behavior of the solution. The reduced number of collocation points, however, provides only a coarse approximation of the underlying PDE and may limit the attainable accuracy. In contrast, optimization on a finer collocation set generally yields a more accurate solution but requires greater computational effort. Motivated by this trade-off, the proposed method gradually increases the number of collocation points while reusing neural network parameters obtained from coarser levels. As illustrated in Figure~\ref{fig:mPINN}, the parameters optimized on a coarse collocation set are transferred to subsequent levels with increasingly dense collocation distributions.

The effectiveness of the proposed strategy relies on the observation that the network parameters obtained at a coarse level typically provide a favorable initialization for the optimization problem at the next finer level. Since these parameters already capture the global behavior of the solution, the subsequent optimization can focus primarily on improving the approximation at newly added collocation points rather than relearning the entire solution from scratch. Consequently, the network is able to resolve finer-scale features more efficiently while avoiding much of the computational cost associated with training from a random initialization. By continuously propagating optimized parameters from coarse to fine levels, the proposed approach accelerates the overall training process. Ultimately, on the finest collocation structure, the method achieves efficient convergence while reducing computational cost and maintaining solution accuracy.
%--------------------------------------

\subsection{Observations on neural network training}

To better understand the rationale behind the proposed coarse-to-fine training strategy, we examine how the number of collocation points influences the optimization problem in PINNs. Let the PDE in (\ref{eq:model}), together with its associated boundary conditions, be written in the abstract form
\begin{equation*}
\cL(U)=\bzero.
\end{equation*}
The corresponding discrete loss function is defined as
\begin{equation}\label{eq:errorLU}
E(\theta)=\sum_{j\in \ti{a set}}|\cL(U_\theta(\bx_j))|^2:=\frac{\lambda_1}{N_r}\sum_{j=1}^{N_r}\left|\mathcal{A}(U_\theta(\bx_j))-F(\bx_j)\right|^2+\frac{\lambda_2}{N_b}\sum_{j=1}^{N_b}
\left|\cB(U_\theta(\bx_j))-G(\bx_j)\right|^2.
\end{equation}
To illustrate the effect of the collocation-point distribution, we compare the following two loss functions:
$$
E_{\ti{coarse}}(\theta)=\sum_{j\in \ti{coarse}}|\cL(U_\theta(\bx_j))|^2\quad\mbox{and}\quad E_{\ti{fine}}(\theta)=\sum_{j\in \ti{fine}}|\cL(U_\theta(\bx_j))|^2.
$$

\begin{figure}
	\centering
	\includegraphics[width=0.90\textwidth]{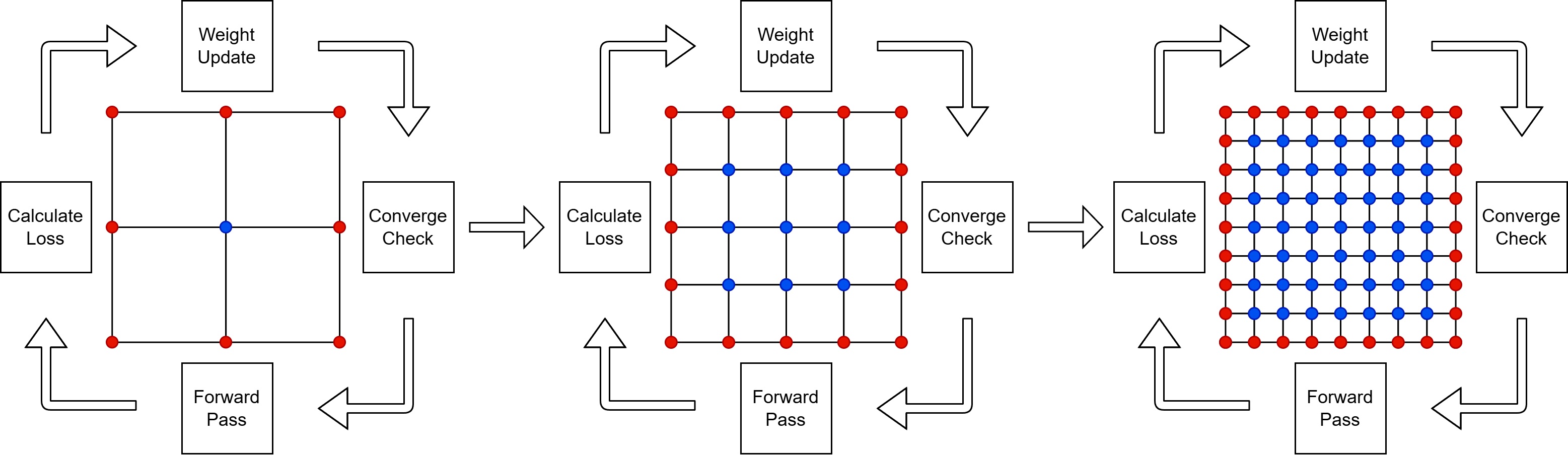}
	\caption{Schematic illustration of the MPU-PINNs procedure, where the parameters learned at a coarse level are used to initialize the training at the next finer level.}
	\label{fig:mPINN}
\end{figure}

It is obvious that using fewer training points reduces the computational cost of each epoch. Based on this observation, the proposed method begins the training process at a coarse level, where the network learns the global behavior of the solution with relatively low computational effort. The parameters trained at this level are then transferred to the subsequent finer level and serve as an effective initial guess. When training process move to a finer level, additional training points are included in the loss function. So, fine level optimization is not identical to the coarse level one, and in principle, the training could be regarded as a new optimization problem. However, in our proposed method, the parameters obtained from the coarse level are used as a good initialization, since they already approximate the global structure of the solution. As a result, the fine level training does not need to learn the entire solution from scratch, but mainly focuses on refining the solution at newly added training points and capturing more detailed local features. This process of transferring parameters from a coarse grid to a fine grid is similar to the strategy of the full multigrid method used in numerical analysis. In such methods, solving the problem on a coarse grid and using it as an initial guess for a finer grid improves the overall efficiency of the solver. When the number of collocation points is small, each training epoch becomes computationally less expensive, since the loss function is evaluated on fewer points. This allows the loss to decrease rapidly at the coarse level with relatively low computational cost. With a smaller number of collocation points, the loss function contains fewer residual terms, making each training epoch computationally cheaper and enabling an efficient reduction of the loss. After obtaining optimal neural network parameters on a coarse set of collocation points, the process can shift to a denser arrangement with a large set of collocation points as shown in Figure \ref{fig:mPINN}. For the transition from a coarse to a finer set of collocation points, two strategies can be employed: a nested method, where the finer grids build upon the coarser structure, or a method where collocation points are randomly selected at each stage of refinement. We call this algorithms as multigrid-based parameter-updated PINNs (MPU-PINNs).

\subsection{Coarse to fine: nested point selection}\label{subsec:pre1}
In the first step of the proposed algorithm, we begin by selecting a small, random set of collocation points within the computational domain $\Omega$. Using PINN, we train the model to satisfy the given PDE (\ref{eq:model}). Once the loss function value falls below a specified tolerance $\epsilon$, additional collocation points are added to the existing set. The PINN is then retrained on the expanded dataset to continue satisfying (\ref{eq:model}). This process is repeated until the number of collocation points reaches the predetermined target number. The figure in Algorithm \ref{algo:pre1} illustrates an example of collocation point sampling process.

\begin{algorithm}[ht!]
\caption{Coarse-to-fine PINN with nested collocation point selection}\label{algo:pre1}
\begin{minipage}[t]{0.7\textwidth}
\begin{algorithmic}[1]
\vspace{10pt}
    \Require A strategy for a sequence of collocation points, an initial collocation set $S$, a target collocation set $T$, residual tolerance $\epsilon$.
    \State $L \Leftarrow 1$
    \While {$S\subsetneq T$}
    \While {$L>\epsilon$}
       \State Train $U_\theta$ on $S$.
       \State Calculate $E(\theta)$ and set $L=E(\theta)$.
    \EndWhile
    \State Generate a set of collocation points, $P\subset T$ with $P\cap S=\emptyset$, to construct a new set $\widetilde{S}$ such that $\widetilde{S}=S\cup P$.
    \State Set $S \Leftarrow \widetilde{S}$ and $L \Leftarrow 1$.
    \EndWhile
    \While {$L>\epsilon$}
        \State Train $U_\theta$ on $S(=T)$.
        \State Calculate $E(\theta)$ and set $L=E(\theta)$.
    \EndWhile
\end{algorithmic}
\end{minipage}
\hfill
\begin{minipage}[t]{0.3\textwidth}
\centering
\vspace{0pt}
\includegraphics[width=0.4\linewidth]{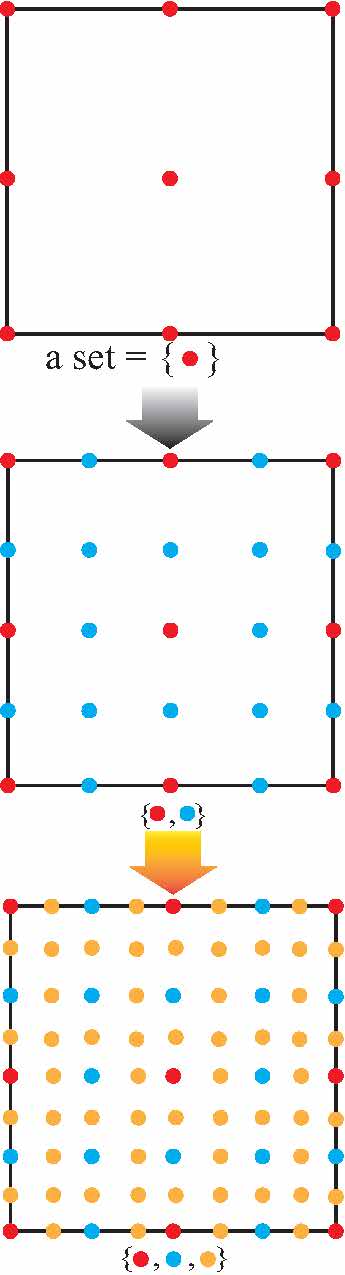}
\end{minipage}
\end{algorithm}

\subsection{Coarse to fine: Random point selection}\label{subsec:pre2}
In this approach, the selection of collocation points follows a purely random process. Starting with a small, randomly chosen set, additional collocation points are sequentially selected until the target number of collocation points is reached. At each stage, the PINN is retrained on the updated dataset, adjusting the neural network parameters, $\theta$, to better approximate the solution of the given PDE as depicted in Algorithm \ref{algo:pre2}. 
\begin{algorithm}[ht!]%\label{algo:2}
\caption{Coarse-to-fine PINN with random collocation point selection}\label{algo:pre2}
\begin{minipage}[t]{0.7\textwidth}
	\begin{algorithmic}[1]
     \vspace{10pt}
		\Require An initial collocation point set $S$, a target collocation point set $T$, residual tolerance $\epsilon$. 
		\State $L \Leftarrow 1$
		\While {$|S| < |T|$, where $|\cdot|$=the number of points in $\cdot$}
		\While {$L>\epsilon$}
		\State  Train $U_\theta$ on $S$.
		\State  Calculate $E(\theta)$ and set $L=E(\theta)$.
		\EndWhile
		\State Generate a new set of collocation points, $\widetilde{S}$, such that $|S| < |\widetilde{S}|$.
		\State Set $S \Leftarrow \widetilde{S}$ and $L \Leftarrow 1$.
		\EndWhile
		\While {$L>\epsilon$}
		\State  Train $U_\theta$ on $S(=T)$.
		\State  Calculate $E(\theta)$ and set $L=E(\theta)$.
		\EndWhile
	\end{algorithmic}
\end{minipage}
\hfill
\begin{minipage}[t]{0.3\textwidth}
\centering
\vspace{0pt}
\includegraphics[width=0.4\linewidth]{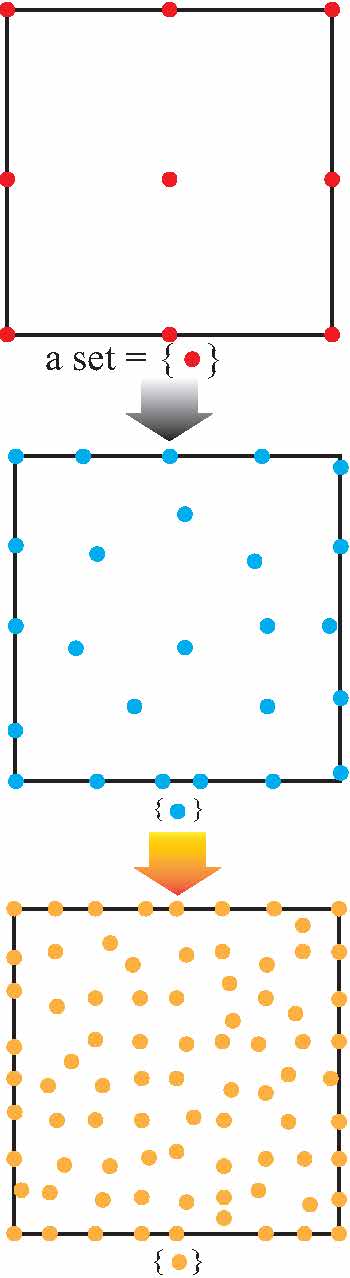}
\end{minipage}
\end{algorithm}

%------------------------------------
\subsection{Efficiency of the coarse-to-fine training strategy}\label{subsec:cost1}

In this section, we discuss why the proposed MPU-PINNs can reduce the training cost compared with conventional PINNs.
Consider a neural network with a fixed architecture and let $P$ denote the number of trainable parameters. Since MPU-PINNs use the same network architecture at every level, the number of trainable parameters remains unchanged throughout the training process. Consequently, the computational savings achieved by MPU-PINNs do not originate from a reduction in model size, but rather from a more efficient optimization strategy.

Let $N_\ell$ denote the number of collocation points at level $\ell$. For a PDE such as the Poisson equation, the residual loss requires the evaluation of second-order derivatives,
\[
R(\bx;\theta)= \mathcal{A}(U_\theta(\bx)) - F(\bx)
\]
at every collocation point. During each optimization step, the following operations are performed:
\begin{enumerate}
\item forward propagation to evaluate $U_\theta$ at all collocation points,
\item automatic differentiation to compute the first- and second-order derivatives appearing in the PDE residual,
\item back propagation to compute the gradient of the loss with respect to all trainable parameters,
\item parameter updates using an optimizer (for example, ADAM, BFGS, Shampoo, and so on).
\end{enumerate}
Among these operations, the dominant cost arises from the repeated evaluation of PDE residuals and the associated automatic differentiation. Since these computations must be carried out at every collocation point, the cost of a single optimization step scales approximately as
$\mathcal{O}(N_\ell P)$, where $N_\ell$ is the number of collocation points and $P$ is the number of trainable parameters. For higher-order PDEs or residuals involving multiple derivatives, the constant hidden in the above estimate becomes larger due to the additional automatic differentiation operations.

For a conventional PINN trained directly on the finest collocation set containing $N_f$ points, the total training cost can be estimated as
\[
{\rm{Cost}}_{\text{\tiny PINN}}\approx K_f\cdot N_L\cdot P,
\]
where $K_f$ denotes the number of optimization iterations required to reach a prescribed accuracy. In contrast, MPU-PINNs employ a sequence collocation sets
with $N_1 < N_2 < \cdots < N_{\mathbb{L}}=N_f$. The total cost becomes
\[
{\rm{Cost}}_{\text{\tiny MPU}} \approx P\sum_{\ell=1}^{\mathbb{L}} K_\ell \cdot N_\ell,
\]
where $K_\ell$ denotes the number of optimization iterations performed at level $\ell$.

At first glance, the multilevel strategy appears to introduce additional computational work because training is performed at multiple levels. However, the crucial observation is that optimization on coarse levels is relatively cheaper since $N_\ell \ll N_f$. More importantly, the parameters obtained from a coarse level provide an effective initialization for the next finer level. As a result, the number of optimization iterations required on fine levels is greatly reduced. In particular, the number of iterations required on the finest level is often much smaller than that of a conventional PINN, $K_{\mathbb{L}} \ll K_f$. Consequently, the expensive optimization process associated with the largest collocation set is substantially shortened. The computational effort is shifted from expensive fine-level iterations to considerably cheaper coarse-level iterations. Since the cost of each iteration grows approximately linearly with the number of collocation points, even a moderate reduction in the number of fine-level iterations can yield a lot of reduction in overall training time.

Another perspective is obtained by viewing the coarse-level training stage as constructing an approximate global representation of the solution. Because the network parameters already capture the dominant low-frequency components of the solution, continued optimization on finer collocation sets primarily refines local features rather than learning the entire solution from scratch. This effect reduces the optimization burden at fine levels and leads to faster convergence. Therefore, although MPU-PINNs do not reduce the number of trainable parameters, they reduce the total computational cost by minimizing the amount of optimization performed on the most expensive collocation sets. This observation is consistent with the numerical experiments presented in the following sections, where substantial reductions in wall-clock training time are observed while maintaining comparable accuracy.

%------------------------------------
\subsection{More evidences of efficiency in the coarse-to-fine training strategy}\label{subsec:cost2}

EXPERIMENTS....

In the following, let's compare with more advanced PINNs...

\begin{itemize}
\item Use examples (test problems) in PIP-IO: reorder the examples 
\item PINN
\item MPU-PINN
\item residual based weighting-PINN
\item residual based weighting-PINN equipped with MPU
\end{itemize}
All using, ADAM + L-BFGS and SOAP.....
\ \\ 

MPU probably has an advantage in using higher-order optimizer due to the small number of derivative for each evaluation

%----------------------------------------
\section{PINN performance comparison with and without neural network parameters updating}\label{sec:pre_examples}

We assess the performance of MPU-PINNs in comparison with conventional PINNs and other PINN-based approaches. More precisely, we consider self-adaptive PINNs (SA-PINNs)~\cite{SA}, which employ adaptive loss weighting, and extended PINNs (XPINNs)~\cite{XPINNs}, which use a domain-decomposition framework. We compare those methods using several PDEs, including two- and three-dimensional Poisson equations, convection-diffusion-reaction equation and Helmholtz equation. These tests aim to show the advantages of incorporating multigrid-based updates into the training process and illustrate how the MPU-PINNs perform in terms of accuracy, convergence speed, and computational efficiency compared to the standard PINNs. For a reliable comparison, all experiments are conducted five times independently, and the average results are reported. We use the relative error, defined as $\|u_{\mathrm{exact}}-u_{\mathrm{approx}}\|/\|u_{\mathrm{exact}}\|$, to evaluate the accuracy of each method.

We consider the following Poisson equation 
\begin{equation}\label{eq:Poisson}
	\left\{\begin{array}{rcll}
		-\Delta u&=& f&\mbox{in}~\Omega,\\
		u&=&g&\mbox{on}~\partial\Omega,\end{array}\right.
\end{equation}
where $\Omega=(0,1)^d,\,d=2,3$, $f$ represents the given source term, and $g$ specifies the prescribed boundary condition on $\partial\Omega$. 
In all Poisson-type equations and convectio--diffusion--reaction equation, a fixed structure of the network that is composed of 5 fully connected hidden layers, each with 50 neurons is used. In addition, the $\tanh$ function is used as the activation function and `Adam' is used as an optimizer with a learning rate of $10^{-3}$. The stopping criterion $\epsilon$ used in Algorithms 1 and 2 is set as $5\times10^{-6}$. The total number of collocation points is set to 263,169. The proposed MPU-PINN starts with 81 collocation points which is gradually increased to 263,169 $(=513\times 513)$ as in Table \ref{tab:CtoF}.

\begin{table}
	\centering
	\caption{Collocation points used in MPU-PINNs from coarse to fine levels}
	\label{tab:CtoF}
	\renewcommand{\arraystretch}{1.1}
	\setlength{\tabcolsep}{10pt}
	\begin{tabular}{cccccccc}
		\toprule
		Level & 1 & 2 & 3 & 4 & 5 & 6 & 7 \\
		\midrule
		$N_r$ in \eqref{eq:errorLU} 
		& 81 & 289 & 1089 & 4225 & 16641 & 66049 & 263169 \\
		$N_b$ in \eqref{eq:errorLU} 
		& 32 & 64 & 128 & 256 & 512 & 1024 & 2048 \\
		\bottomrule
	\end{tabular}
\end{table}

\begin{figure}
	\centering
	\includegraphics[width=0.9\textwidth]{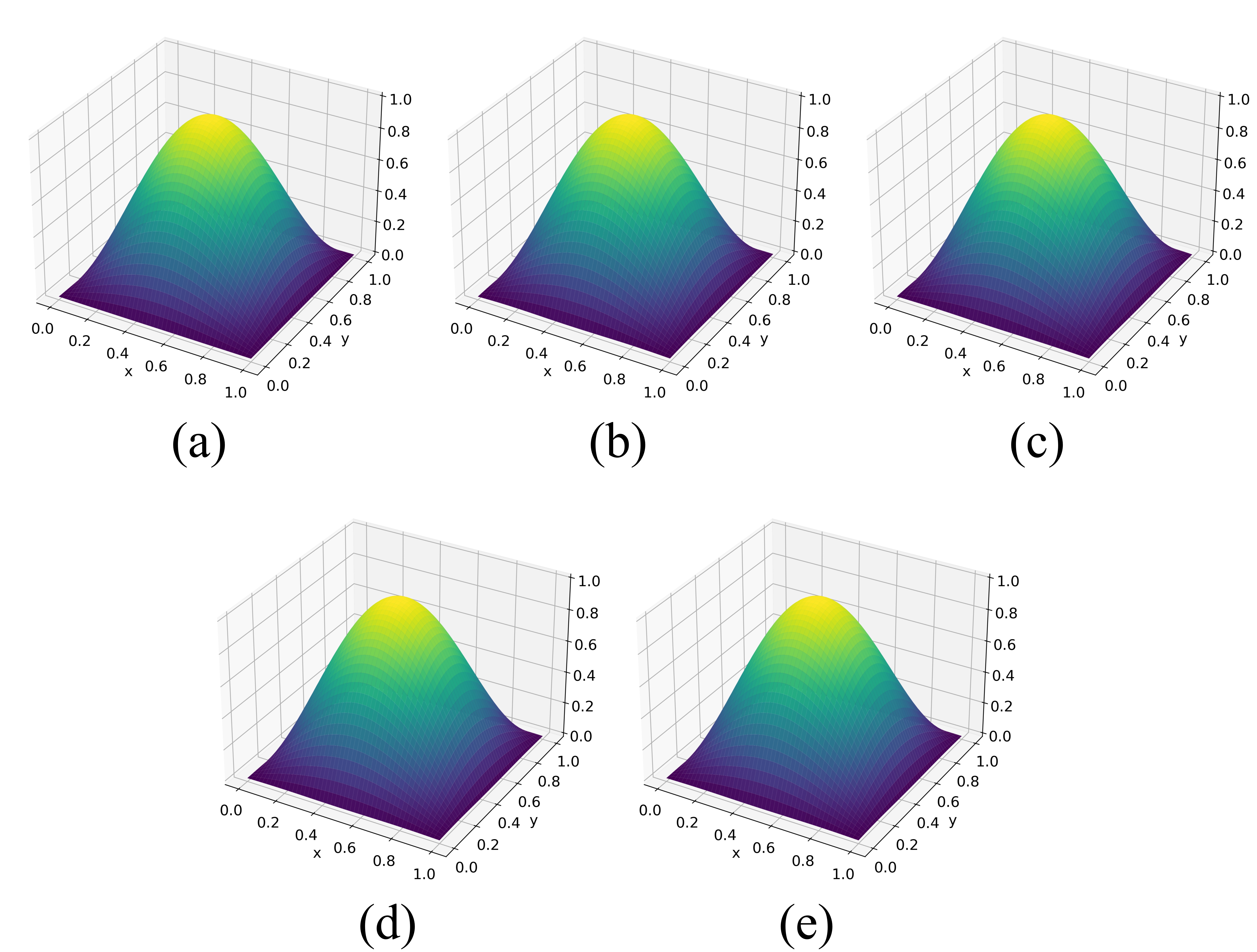}
	\caption{Exact and predicted solutions for Example~\ref{exam:ex1}: (a) exact solution (b) PINNs (c) SA-PINNs (d) MPU-PINN Algorithm 1 (e) MPU-PINN Algorithm 2}
	\label{fig:ex1}
\end{figure}

\begin{example}\label{exam:ex1}
	We set $\Omega=(0,1)^2$ and the exact $u(x,y)=\sin(\pi x)\sin(\pi y)$ that yield $f(x,y)=2\pi^2 \sin(\pi x)\sin(\pi y)$ in $\Omega$ and $g(x,y)=0$ on $\Gamma$. The weight coefficients in (\ref{eq:errorLU}) are set to $\lambda_1=1$ and $\lambda_2=10^{-2}$.
\end{example}
Table \ref{table:ex1_timeNerr} shows that MPU-PINN outperforms conventional PINN while achieve a similar level of accuracy. The time required for the MPU-PINN algorithms to reduce the loss function below a given tolerance is nearly 5$\%$ shorter compared to the conventional PINN with the relative error remains comparable across all tests. Table \ref{table:ex1_epoch} presents the total number of epochs required for MPU-PINN to reduce the loss function below the tolerance $\epsilon$ at each level whereas conventional PINN requires 8344 epochs. From Table \ref{table:ex1_epoch}, it is clear that when the solution is smooth, the neural network parameters are well-optimized on coarse grids which provide good initial estimates for the subsequent finer grids. This is indicated by the decreasing number of epochs required as the grid becomes progressively finer. Figure \ref{fig:ex1} exhibits the exact solution alongside the approximations obtained from various PINN algorithms. It can be observed that all methods successfully capture the exact solution.

\begin{table}[h]
	\centering
	\caption{Comparison} \label{table:ex1_timeNerr}
	\renewcommand{\arraystretch}{1.15}
	\setlength\tabcolsep{12pt}
	\begin{tabular}{llcc}
		\toprule
		Category & Method & Relative error & Elapse time (seconds) \\
		\midrule
		\multirow{2}{*}{Baseline}
		& Conventional PINN & 3.8101e-03 & 1650.9765 \\
		& SA-PINN           & 7.4560e-03 & 784.9539  \\
		\midrule
		\multirow{2}{*}{Proposed}
		& MPU-PINN Algorithm 1   & 3.2818e-03 & 51.5128 \\
		& MPU-PINN Algorithm 2   & 3.0146e-03 & 67.7841  \\
		\bottomrule
	\end{tabular}
\end{table}

\begin{table}[h]
	\centering
	\caption{Epoch distribution over training levels in Example~\ref{exam:ex1}.}
	\label{table:ex1_epoch}
	\renewcommand{\arraystretch}{1.15}
	\setlength{\tabcolsep}{7pt}
	\begin{tabular}{llccccccc}
		\toprule
		\multirow{2}{*}{Category} 
		& \multirow{2}{*}{Method} 
		& \multicolumn{7}{c}{Training level} \\
		\cmidrule(lr){3-9}
		& & 1 & 2 & 3 & 4 & 5 & 6 & 7 \\
		\midrule
		\multirow{2}{*}{Baseline}
		& Conventional PINNs & -- & -- & -- & -- & -- & -- & 8344.4 \\
		& SA-PINNs           & -- & -- & -- & -- & -- & -- & 3856 \\
		\midrule
		\multirow{2}{*}{Proposed}
		& MPU-PINNs Algorithm 1 & 5856.4 & 1231.0 & 1 & 1 & 1 & 1 & 1 \\
		& MPU-PINNs Algorithm 2 & 6790.8 & 1963.4 & 116.6 & 6.6 & 1 & 3.4 & 1.5 \\
		\bottomrule
	\end{tabular}
\end{table}

%\begin{table}
%\caption{total number of epochs to achieve loss reduction below tolerance$\epsilon$}\label{table:ex1_epoch}
%\renewcommand{\arraystretch}{1.1}
%\setlength\tabcolsep{9pt}
%\begin{tabular}{lccccccc}
%\toprule
%method&1&2&3&4&5&6&7\\ \midrule
%Conventional PINNs		   & - & - & - & - & - & - & 8344.4 \\
%SA-PINNs		           & - & - & - & - & - & - & 3856 \\
%MPU-PINNs Algorithm 1      & 5856.4 & 1231.0 & 1 & 1 & 1 & 1 & 1 \\
%MPU-PINNs Algorithm 2      & 6790.8 & 1963.4 & 116.6 & 6.6 & 1 & 3.4 & 1.5 \\
%\botrule
%\end{tabular}
%\end{table}

Now, we test our MPU-PINN on following convection--diffusion--reaction problem. 
\begin{example}\label{exam:ex2}
	\begin{equation}
		\begin{cases}
			-\epsilon \Delta u + \mathbf{b}\cdot \nabla u + \sigma u = f,
			& \quad \text{in} \quad \Omega, \\
			\quad \,\, u =  u_{exact}, 
			& \quad \text{on} \quad \partial\Omega.
		\end{cases}
	\end{equation}
\end{example}

The parameters are chosen as $\epsilon=10^{-8}$, $\mathbf{b}=(2,3)^T$, and $\sigma=2$, which makes the problem strongly convection-dominated. 
It yields exact solution $u_{\mathrm{exact}}(x,y) = 16x(1-x)y(1-y) \left( \frac{1}{2} + \frac{\arctan\left(200\left(r_0^2-(x-x_c)^2-(y-y_c)^2\right)\right)}{\pi} \right)$. the weight coefficients are set as $(\lambda_1,\lambda_2)=(1,10^{-3})$. Unlike previous example, it has steep transition  training is started from level 4 with $N_r=4225$, rather than from the level 1. This is because the exact solution contains a sharp internal transition near the center of domain, and excessively coarse training points, like level 1, may fail to represent this localized structure. Therefore, we think a sufficiently refined points are required so that the network can capture the steep transition region from the early stage of training. Similar to the results in Example \ref{exam:ex1}, Table \ref{table:ex2_timeNerr} depicts the robust performance of MPU-PINN. Moreover, Table \ref{table:ex2_epoch} indicates that finer grids require fewer training epochs that results in shorter overall computational time (conventional PINN takes 31848 epochs). Generally, having a larger number of collocation points rises the complexity of the neural network training which would typically require more training epochs to achieve convergence. However, with MPU-PINN, this challenge is mitigated because the network's parameters have already been effectively optimized in earlier stages. The optimized parameters at these earlier stages serve as a strong initial guess for the finer grids, effectively reducing the search space for the neural network during subsequent training. This concept is aligned with transfer learning in neural networks, where knowledge gained from simpler tasks is used to accelerate the training of more complex tasks. As a result, even though the number of collocation points increases, the neural network requires fewer epochs to bring the loss function below the tolerance level $\epsilon$ because it starts with parameters already close to the optimal solution. This efficient transfer of information from earlier levels accelerates the convergence process, as reflected in the reduced number of epochs seen in Table \ref{table:ex2_epoch}. 

\begin{table}[h]
	\centering
	\caption{Comparison} \label{table:ex2_timeNerr}
	\renewcommand{\arraystretch}{1.15}
	\setlength\tabcolsep{12pt}
	\begin{tabular}{llcc}
		\toprule
		Category & Method & Relative error & Elapse time (seconds) \\
		\midrule
		\multirow{2}{*}{Baseline}
		& Conventional PINN & 6.8724e-03 & 6560.9304 \\
		& SA-PINN           & 8.6071e-03 & 4361.6970 \\
		\midrule
		\multirow{2}{*}{Proposed}
		& MPU-PINN Algorithm 1   & 6.8151e-03 & 502.2478  \\
		& MPU-PINN Algorithm 2   & 7.5682e-03 & 636.9776  \\
		\bottomrule
	\end{tabular}
\end{table}

\begin{table}[h]
	\centering
	\caption{Epoch distribution over training levels in Example~\ref{exam:ex2}.}
	\label{table:ex2_epoch}
	\renewcommand{\arraystretch}{1.15}
	\setlength{\tabcolsep}{7pt}
	\begin{tabular}{ll *{4}{>{\centering\arraybackslash}p{1.3cm}}}
		\toprule
		\multirow{2}{*}{Category} 
		& \multirow{2}{*}{Method} 
		& \multicolumn{4}{c}{Training level} \\
		\cmidrule(lr){3-6}
		& & 4 & 5 & 6 & 7 \\
		\midrule
		\multirow{2}{*}{Baseline}
		& Conventional PINNs & - & - & - & 31848.4 \\
		& SA-PINNs           & - & - & - & 21207.2 \\
		\midrule
		\multirow{2}{*}{Proposed}
		& MPU-PINNs Algorithm 1 & 29695.4 & 9081.4 & 1087.4 & 55.8 \\
		& MPU-PINNs Algorithm 2 & 31786.6 & 10512.6 & 2588.6 & 494.8 \\
		\bottomrule
	\end{tabular}
\end{table}

%\begin{table}
%	\caption{total number of epochs to achieve loss reduction below tolerance $\epsilon$}\label{table:ex2_epoch}
%	\renewcommand{\arraystretch}{1.1}
%	\setlength\tabcolsep{9pt}
%	\begin{tabular}{lcccc}
	%		\toprule
	%		level&4&5&6&7\\ \midrule
	%		Conventional PINNs		   & - & - & - & 31848.4 \\
	%		SA-PINNs		           & - & - & - & 21207.2 \\
	%		MPU-PINNs Algorithm 1      & 29695.4 & 9081.4 & 1087.4 & 55.8 \\
	%		MPU-PINNs Algorithm 2      & 31786.6 & 10512.6 & 2588.6 & 494.8 \\
	%		\botrule
	%	\end{tabular}
%\end{table}

\begin{figure}
	\centering
	\includegraphics[width=0.9\textwidth]{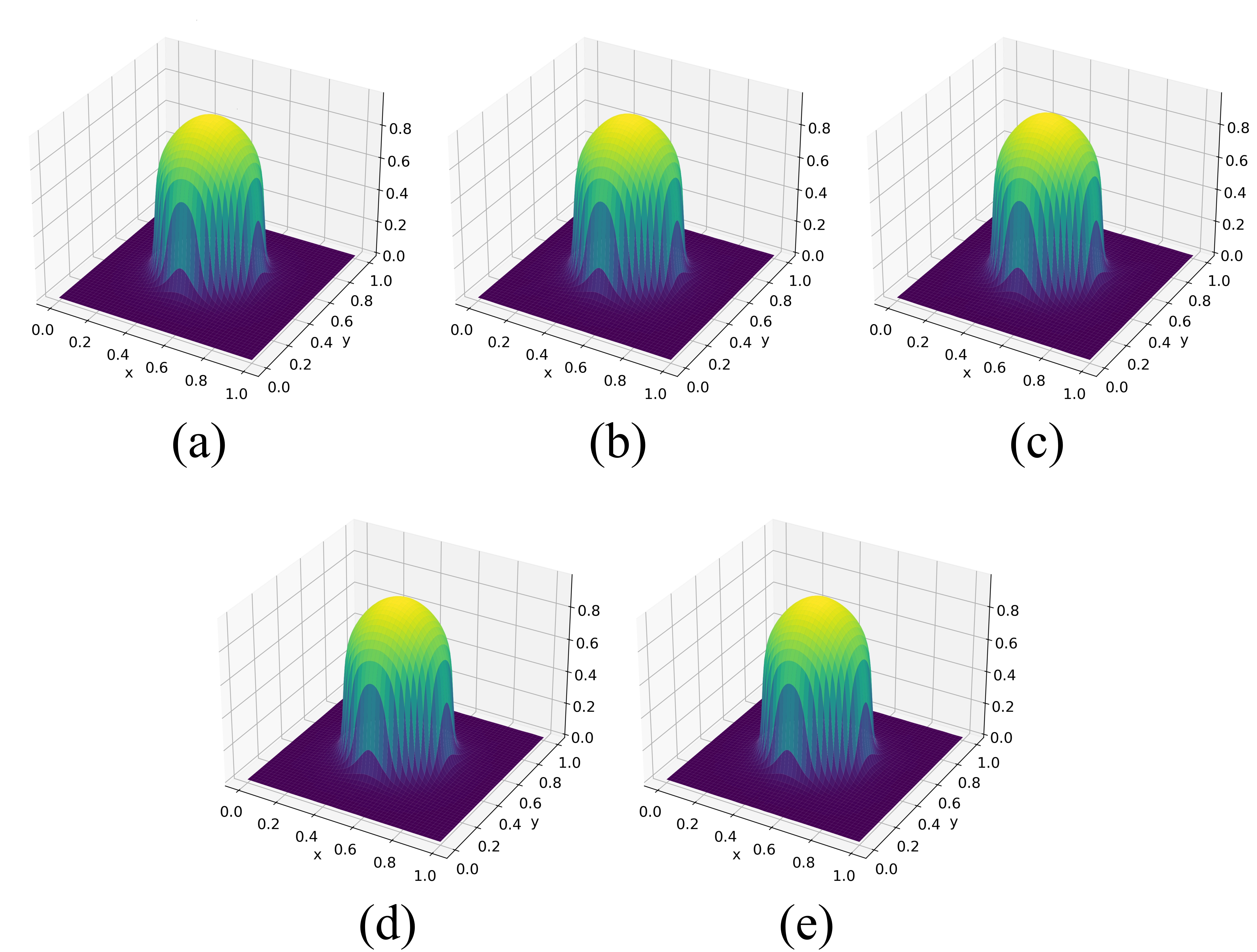}
	\caption{Exact and predicted solutions for Example~\ref{exam:ex2}: (a) exact solution (b) PINNs (c) SA-PINNs (d) MPU-PINN Algorithm 1 (e) MPU-PINN Algorithm 2}
	\label{fig:exa2}
\end{figure}

Figure \ref{fig:exa2} also displays all PINN algorithms provides good approximations. 

\section{Application to High-Frequency Problems}\label{sec:high}

As is common with traditional numerical methods, capturing high-frequency components necessitates the use of a highly refined mesh. This requirement increases the computational cost. From this perspective, high-frequency problems provide a suitable test case for evaluating the effectiveness of MPU-PINNs. Since MPU-PINNs are designed to accelerate convergence on fine-level through coarse-to-fine training, they are expected to be useful methods. In addition to the increased computational cost, neural networks has limitation, where low-frequency components are learned first than high-frequency components, which called spectral bias. This phenomenon make neural networks difficult approximate oscillatory solutions. To reduce limitation, we adopt sinusoidal activation strategy, which use sinusoidal function at the first layer and sin activation functions are used in other layers. This strategy is expected to improve the trainability of input gradient of neural network at the beginning of optimization, thereby reducing the possibility of being trapped in undesirable local minima.

\begin{figure}
	\centering
	\includegraphics[width=0.9\textwidth]{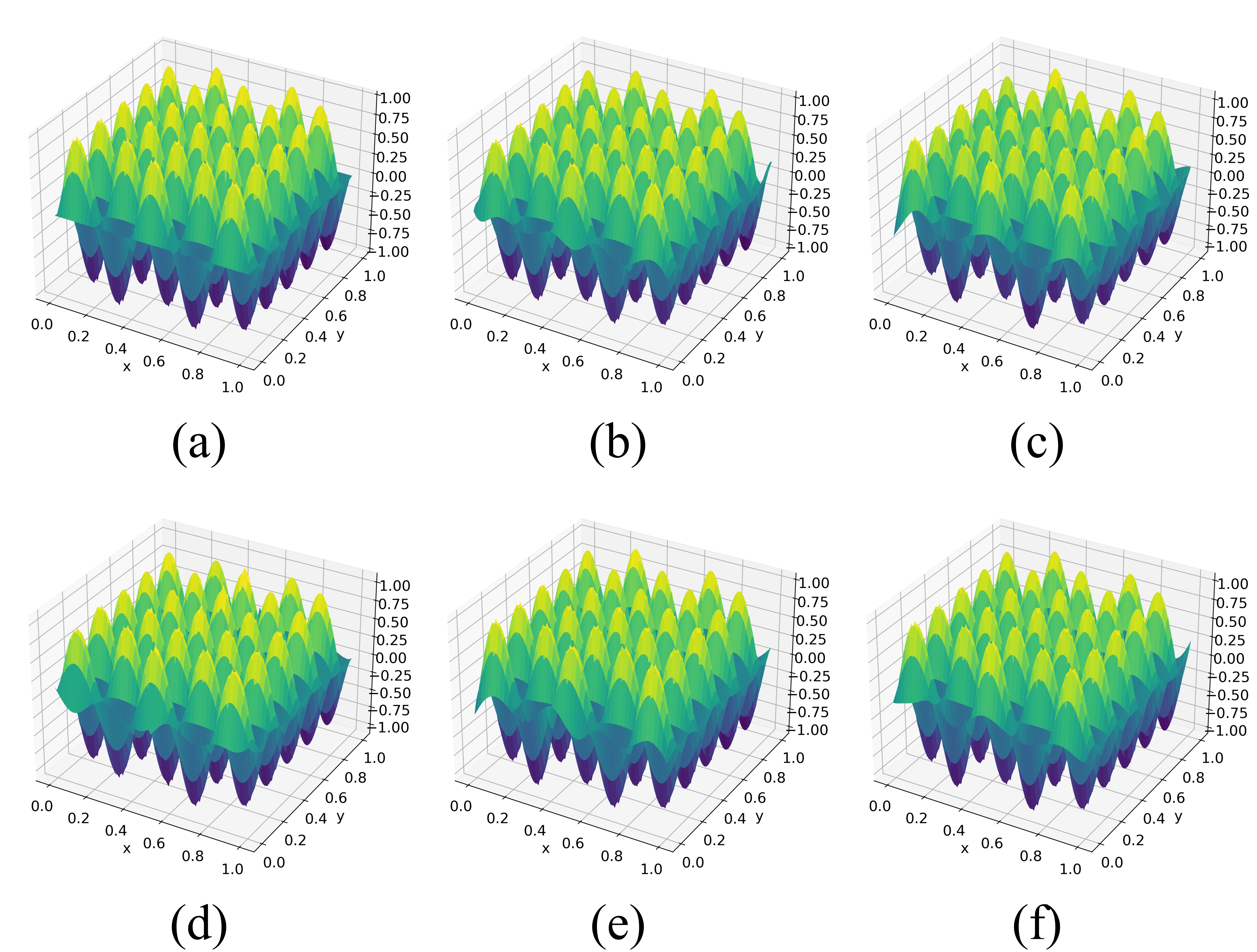}
	\caption{Exact and predicted solutions for Example~\ref{exam:ex3}: (a) exact solution (b) PINNs (c) SA-PINNs (d) XPINNs (e) MPU-PINN Algorithm 1 (f) MPU-PINN Algorithm 2}
	\label{fig:ex3}
\end{figure}

\begin{example}\label{exam:ex3}
	In this example, we revisit \ref{eq:Poisson} and set $\Omega=(0,1)^2$ with the exact solution $u(x,y)=\sin(5\pi x)\sin(10\pi y)$ which leads to $f(x,y)=125\pi^2 \sin(5\pi x)\sin(10\pi y)$ and boundary condition $g(x,y)=0$.
\end{example}
While the solution is smooth, it exhibits relatively high frequencies compared to Example \ref{exam:ex1}. Therefore, starting with level 1, which has few training points ($Nr=81$), is leff efficient in this case. From this observation, we start training from level 4. Conventional PINNs and XPINNs reached the maximum of 100,000 epochs without reducing the loss function below the specified tolerance, whereas SA-PINNs required 71959 epochs. Similar to previous examples, required number of epochs for each level is gradually decrease. This suggests that the parameter transfer between levels is effective even for problems with high frequency. In Table ~\ref{table:ex3_timeNerr}, we can observe MPU-PINNs significantly reduce the computational cost compared with the baseline methods. MPU-PINNs uses less than 7\% of the elapsed time required by conventional PINNs and only about 9--11\% of that required by SA-PINNs and XPINNs, while achieving comparable relative errors. Figure~\ref{fig:ex3} shows that all methods reproduce the overall profile of the exact solution over the entire domain. But they does not satisfy the boundary condition exactly. This behavior is reflected in the relative errors reported in Table~\ref{table:ex3_timeNerr}, which remain on the order of $10^{-1}$. This observation motivates the improved strategy introduced in the following section.

\begin{table}[h]
	\centering
	\caption{Epoch distribution over training levels in Example~\ref{exam:ex3}.}
	\label{table:ex3_epoch}
	\renewcommand{\arraystretch}{1.15}
	\setlength{\tabcolsep}{7pt}
	\begin{tabular}{ll *{4}{>{\centering\arraybackslash}p{1.3cm}}}
		\toprule
		\multirow{2}{*}{Category} 
		& \multirow{2}{*}{Method} 
		& \multicolumn{4}{c}{Training level} \\
		\cmidrule(lr){3-6}
		& & 4 & 5 & 6 & 7 \\
		\midrule
		\multirow{2}{*}{Baseline}
		& Conventional PINNs & - & - & - & 100000 \\
		& SA-PINNs           & - & - & - & 71959.6 \\
		& XPINNs		           	   & - & - & - & 100000 \\
		\midrule
		\multirow{2}{*}{Proposed}
		& MPU-PINNs Algorithm 1 & 99\,099 & 33372.8 & 1 & 85 \\
		& MPU-PINNs Algorithm 2 & 100000 & 48488 & 460.6 & 1 \\
		\bottomrule
	\end{tabular}
\end{table}

\begin{table}[h]
	\centering
	\caption{Comparison} \label{table:ex3_timeNerr}
	\renewcommand{\arraystretch}{1.15}
	\setlength\tabcolsep{12pt}
	\begin{tabular}{llcc}
		\toprule
		Category & Method & Relative error & Elapse time (seconds) \\
		\midrule
		\multirow{2}{*}{Baseline}
		& Conventional PINN & 2.5716e-01 & 21131.3467 \\
		& SA-PINN           & 2.5716e-01 & 11826.8132 \\
		& XPINNs           & 3.6227e-01 & 12232.3401 \\
		\midrule
		\multirow{2}{*}{Proposed}
		& MPU-PINN Algorithm 1   & 2.3435e-01 & 1037.3344  \\
		& MPU-PINN Algorithm 2   & 2.5014e-01 & 1267.9327  \\
		\bottomrule
	\end{tabular}
\end{table}

\subsection{A new strategy of scaling in loss function}\label{sec:scaling}

Solving problems such as the Helmholtz equation, which depends on the wave number, poses considerable obstacles even for conventional numerical methods. Similarly, the vanilla PINN algorithms also face difficulties in these cases. In this section, we propose introducing a the loss function which apply scaling to mitigate these challenges and obtain better approximate solutions using PINNs.

We recall the Poisson equation introduced in Example \ref{exam:ex3} and rewrite it as follows:
\begin{equation*}\label{eq:ex3_rewrite}
	\left\{\begin{array}{rcll}
		-\Delta u&=& k^2 \tilde f&~~~\mbox{in}~\Omega,\\
		u&=&0&~~~\mbox{on}~\partial\Omega,\end{array}\right.
\end{equation*}
where $\tilde f(x,y) =\sin(5\pi x)\sin(10\pi y)$ and $k=5\sqrt{5}\pi$. 
This can be reformulated as an equivalent system of equations:
\begin{equation}\label{eq:ex3_scale}
	-\frac{1}{k^2}\Delta u= \tilde f, ~~~\mbox{in}~\Omega
\end{equation}
with $u=0$ on $\partial\Omega$. Then we define the discrete loss function as 
\begin{align}\label{eq:ex3_loss}
	E(\theta)=&\frac{1}{N_r}\sum_{j=1}^{N_r}\lambda_{eq}\left|-\frac{1}{k^2}\Delta u - \tilde f\right|^2 \,\, + \,\, \frac{1}{N_b}\sum_{j=1}^{N_b}\lambda_{bc}\left|u(\bx_j)\right|^2,
\end{align}
where $\bx_j=(x_j,y_j)$.
The another settings for MPU-PINNs remain identical to \ref{exam:ex3}, except for loss formulation. We assign same weight coefficients as follows: $\lambda_{eq}=0.1$, and $\lambda_{bc}=1.0$. Tables~\ref{table:ex3_scale_timeNerr} and~\ref{table:ex3_scale_epoch} shows effectiveness of the scaling strategy. Compared with without scaling formulation, scaling equations and inducing an appropriate corresponding error function can enhance the overall performance of the PINN. Nevertheless, the results also show that sufficiently refined training points are required as the oscillatory behavior becomes more pronounced. Therefore, the initial training level should be carefully selected for problems with large wave numbers. To examine this issue more systematically, we consider the Helmholtz equation in the following section.

\begin{table}[h]
	\centering
	\caption{Effect of scaling on the elapsed time and relative error of MPU-PINNs for Example~\ref{exam:ex3}.}
	\label{table:ex3_scale_timeNerr}
	\renewcommand{\arraystretch}{1.2}
	\setlength{\tabcolsep}{9pt}
	\begin{tabular}{l c c c c}
		\toprule
		Method & Sampling & Scaling & Elapsed time (s) & Relative error \\
		\midrule
		\multirow{4}{*}{MPU-PINNs}
		& \multirow{2}{*}{Uniform}
		& Without scaling & 1037.3344 & 2.3435e-01 \\
		& & With scaling    & 417.8385  & 2.2195e-02 \\
		\cmidrule(lr){2-5}
		& \multirow{2}{*}{Random}
		& Without scaling & 1267.9327 & 2.5014e-01 \\
		& & With scaling    & 420.1623  & 2.9938e-02 \\
		\bottomrule
	\end{tabular}
\end{table}

\begin{table}[h]
	\centering
	\caption{Epoch distribution across training levels for MPU-PINNs with and without scaling in Example~\ref{exam:ex3}.}
	\label{table:ex3_scale_epoch}
	\renewcommand{\arraystretch}{1.2}
	\setlength{\tabcolsep}{9pt}
	\begin{tabular}{l c c c c c c}
		\toprule
		Method & Sampling & Scaling & Level 4 & Level 5 & Level 6 & Level 7 \\
		\midrule
		\multirow{4}{*}{MPU-PINNs}
		& \multirow{2}{*}{Uniform}
		& Without scaling & 99099  & 33372.8 & 1     & 85 \\
		& 
		& With scaling    & 48882.6 & 1.6     & 1     & 1 \\
		\cmidrule(lr){2-7}
		& \multirow{2}{*}{Random}
		& Without scaling & 100000  & 48488   & 460.6 & 1 \\
		& 
		& With scaling    & 49244.2 & 3145.4  & 2.6   & 16.6 \\
		\bottomrule
	\end{tabular}
\end{table}

%-----------------------------------------------------

\subsection{Helmholtz equation}\label{sec:Helmholtz}

The Helmholtz equation presented below models wave propagation and resonance phenomena in physical systems with an external source. It describes situations involving forced oscillations or wave interactions with sources in a bounded domain:

\begin{example}\label{exam:ex4}
	\begin{equation}\label{eq:Helmholtz}
		\begin{cases}
			-\Delta u - k^2 u = k^2 \tilde{f},
			& \quad \text{in} \quad \Omega, \\
			\quad \,\, u = 0,
			& \quad \text{on} \quad \partial\Omega.
		\end{cases}
	\end{equation}
\end{example}

Here, $k$ denotes the wave number that is related to the wavelength by $\lambda=2\pi/k$. We test the conventional PINN, SA-PINNs, XPINNs and proposed MPU-PINNs that starts with mesh $1/32$ and moves to $1/256$ with $k=8\pi$ and $\tilde f = \sin(8\pi x)\sin(8\pi y)$. Similar to Example~\ref{exam:ex3}, we also scale given equation. The corresponding loss function is defined as
\begin{equation}\label{eq:Helm_ori}
	E(\theta)=\lambda_{eq}\cdot\frac{1}{N_r}\sum_{j=1}^{N_r}\left|- \frac{1}{k^2} \Delta u(\bx_j)- u(\bx_j)-\tilde f(\bx_j)\right|^2+\lambda_{bc}\cdot\frac{1}{N_b}\sum_{j=1}^{N_b}\left|u(\bx_j)\right|^2,
\end{equation}
where $\lambda_{ed} = 1.0$, $\lambda_{bc} = 0.01$. We use a neural network with five hidden layers and 100 neurons per layer. Training is terminated once the loss value reaches $5\times 10^{-7}$. In XPINNs, we use two subdomains, and each subdomain is assigned a local network with four hidden layers and 81 neurons per layer so that the total number of trainable parameters is comparable to that of conventional PINNs.

Table~\ref{table:ex4_epoch} illustrates the training behavior of each method. Conventional PINNs and SA-PINNs require a large number of epochs at the finest level, and XPINNs reaches the maximum epoch. On the other hand, both MPU-PINNs algorithms satisfy the stopping criterion at all training levels before reaching the maximum number of epochs. This confirms that the proposed multilevel training strategy provides a more efficient optimization process for this high-frequency problem. As shown in Table~\ref{table:ex4_timeNerr}, the proposed MPU-PINNs achieve relative errors that are close to those of conventional PINNs, while requiring much shorter computational time. In particular, MPU-PINNs Algorithm 1 and Algorithm 2 reduce the elapsed time from $10323.7600$ seconds for conventional PINNs to $1014.3697$ and $885.3785$ seconds, respectively. This indicates that our coarse-to-fine training strategy can substantially improve computational efficiency without a significant loss of accuracy. Also compared with SA-PINNs and XPINNs, MPU-PINNs also show a clear advantage in terms of elapsed time. Although SA-PINNs and XPINNs require less computational time than conventional PINNs, their relative errors are larger in this example. These results suggest that MPU-PINNs are well suited for the Helmholtz equation, where the oscillatory solution structure makes the training process more challenging.

\begin{figure}
	\centering
	\includegraphics[width=0.9\textwidth]{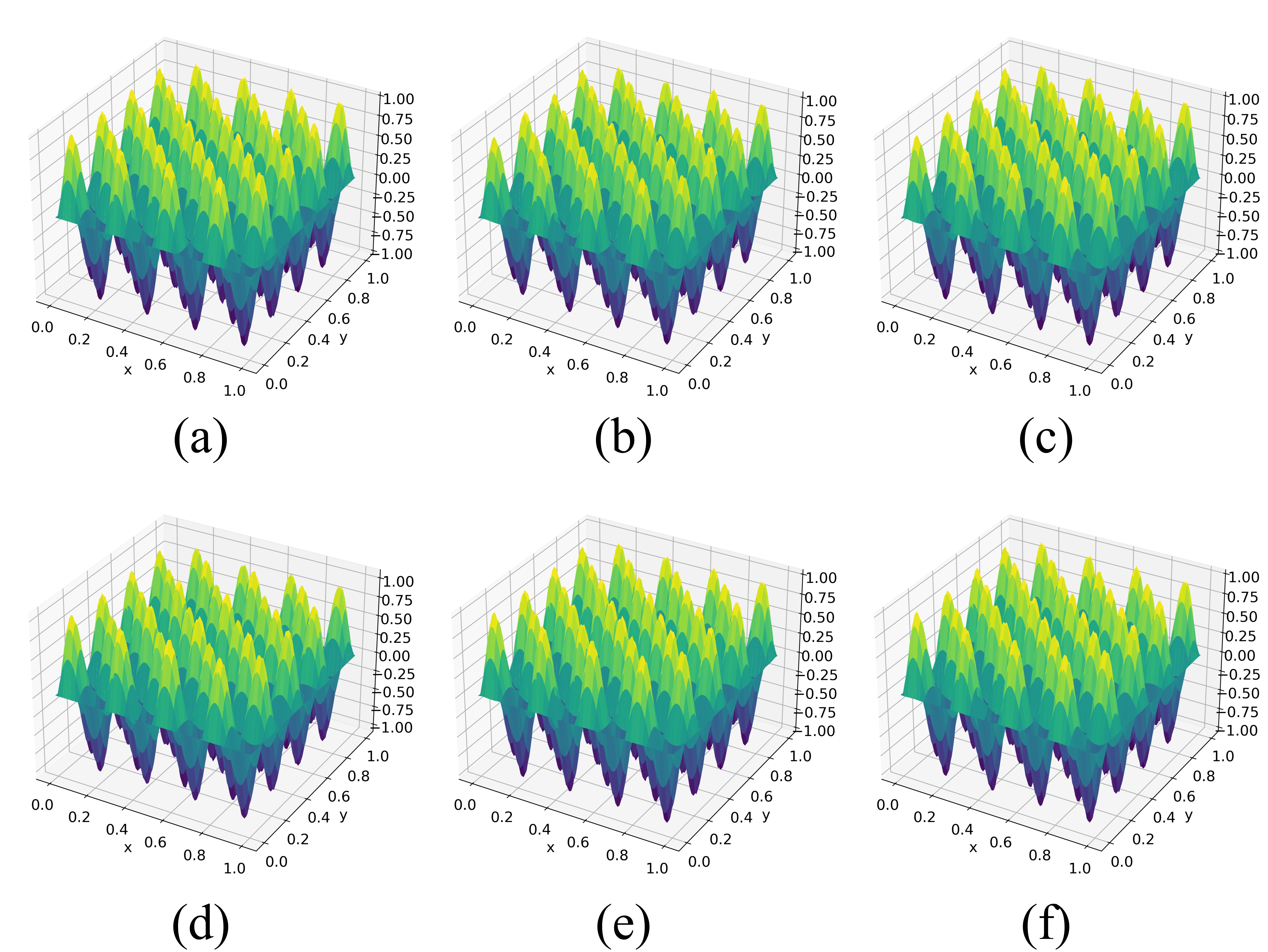}
	\caption{Exact and predicted solutions for Example~\ref{exam:ex4}: (a) exact solution (b) PINNs (c) SA-PINNs (d) XPINNs (e) MPU-PINN Algorithm 1 (f) MPU-PINN Algorithm 2}
	\label{fig:ex4}
\end{figure}

\begin{table}[h]
	\centering
	\caption{Epoch distribution over training levels in Example~\ref{exam:ex4}.}
	\label{table:ex4_epoch}
	\renewcommand{\arraystretch}{1.15}
	\setlength{\tabcolsep}{7pt}
	\begin{tabular}{ll *{4}{>{\centering\arraybackslash}p{1.3cm}}}
		\toprule
		\multirow{2}{*}{Category} 
		& \multirow{2}{*}{Method} 
		& \multicolumn{4}{c}{Training level} \\
		\cmidrule(lr){3-6}
		& & 3 & 4 & 5 & 6 \\
		\midrule
		\multirow{2}{*}{Baseline}
		& Conventional PINNs & -- & -- & -- & 97730 \\
		& SA-PINNs           & -- & -- & -- & 77253 \\
		& XPINNs		           	   & -- & -- & -- & 100000 \\
		\midrule
		\multirow{2}{*}{Proposed}
		& MPU-PINNs Algorithm 1 & 62386.6 & 37850.4 & 4817.8 & 1028.6 \\
		& MPU-PINNs Algorithm 2 & 61464 & 42678.2 & 3582.8 & 267 \\
		\bottomrule
	\end{tabular}
\end{table}

\begin{table}[h]
	\centering
	\caption{Comparison} \label{table:ex4_timeNerr}
	\renewcommand{\arraystretch}{1.15}
	\setlength\tabcolsep{12pt}
	\begin{tabular}{llcc}
		\toprule
		Category & Method & Relative error & Elapse time (seconds) \\
		\midrule
		\multirow{2}{*}{Baseline}
		& Conventional PINN & 1.8050e-02  & 10323.7600 \\
		& SA-PINN           & 4.4598e-02 & 8247.2385 \\
		& XPINNs           & 8.4218e-02 & 7890.7014 \\
		\midrule
		\multirow{2}{*}{Proposed}
		& MPU-PINN Algorithm 1   & 2.3069e-02  & 1014.3697  \\
		& MPU-PINN Algorithm 2   & 2.1243e-02 & 885.3785  \\
		\bottomrule
	\end{tabular}
\end{table}

\newpage

\section{Conclusion}\label{sec:conclude}
In this paper, we proposed multigrid-based parameter-updated PINNs (MPU-PINNs) to improve the training efficiency of physics-informed neural networks. 
The main idea of the proposed method is to train a single neural network in a coarse-to-fine manner by progressively increasing the number of collocation points. 
At each level, the parameters obtained from the previous level are transferred to the next finer level as an effective initialization. 
This allows the network to first capture the global structure of the solution with a relatively small number of training points and then refine the solution using additional points at finer levels.

Several benchmark problems were considered to evaluate the performance of the proposed method, including two- and three-dimensional Poisson equations, a Helmholtz equation, and a convection--diffusion--reaction equation. 
The numerical results showed that MPU-PINNs significantly reduce the computational time compared with conventional PINNs and representative PINN-based methods such as SA-PINNs and XPINNs, while maintaining comparable accuracy. 
In particular, the epoch distributions demonstrated that the proposed coarse-to-fine strategy reduces the training burden at finer levels and avoids training the network from random initialization at each refinement stage.

For high-frequency problems, we also introduced a scaling strategy to improve the training behavior of MPU-PINNs. 
The results showed that the scaled formulation reduces the relative error and accelerates convergence compared with the unscaled formulation. 
This indicates that an appropriate reformulation of the governing equation and the corresponding loss function can play an important role in improving the performance of PINNs for oscillatory problems.

The convection--diffusion--reaction example further showed that the choice of the initial training level is important when the solution contains a localized steep transition. 
In such cases, starting from an excessively coarse level may not provide sufficient information to resolve the sharp transition region. 
Therefore, a moderately refined initial level can be more effective for layer-type or convection-dominant problems.

Overall, the proposed MPU-PINNs provide an efficient single-network training framework for solving PDEs with PINNs. 
By combining parameter transfer with progressive refinement of collocation points, the method improves computational efficiency while preserving the simplicity of the standard PINN architecture. 
Future work will consider a wider range of PDE problems to further examine the general applicability of the MPU-PINNs. In addition, theoretical analysis will be pursued to better understand the convergence properties and the effect of the coarse-to-fine parameter transfer strategy.

\printcredits

\bibliographystyle{cas-model2-names}
%\bibliographystyle{unsrtnat}

\bibliography{cas-refs}

@article{Raissi,
  author  = {Raissi, Maziar and Perdikaris, Paris and Karniadakis, George Em},
  title   = {Physics-informed neural networks: A deep learning framework for solving forward and inverse problems involving nonlinear partial differential equations},
  journal = {Journal of Computational Physics},
  volume  = {378},
  pages   = {686--707},
  year    = {2019},
  doi     = {10.1016/j.jcp.2018.10.045}
}

@article{FF,
  author  = {Wang, Sifan and Wang, Hanwen and Perdikaris, Paris},
  title   = {On the eigenvector bias of Fourier feature networks: From regression to solving multi-scale {PDEs} with physics-informed neural networks},
  journal = {Computer Methods in Applied Mechanics and Engineering},
  volume  = {384},
  pages   = {113938},
  year    = {2021},
  doi     = {10.1016/j.cma.2021.113938}
}

@inproceedings{Instable_PINN,
  author    = {Krishnapriyan, Aditi S. and Gholami, Amir and Zhe, Shandian and Kirby, Robert and Mahoney, Michael W.},
  title     = {Characterizing possible failure modes in physics-informed neural networks},
  booktitle = {Advances in Neural Information Processing Systems},
  volume    = {34},
  pages     = {26548--26560},
  year      = {2021},
  url       = {https://proceedings.neurips.cc/paper/2021/hash/df438e5206f31600e6ae4af72f2725f1-Abstract.html}
}

@article{NTK,
  author  = {Wang, Sifan and Yu, Xinling and Perdikaris, Paris},
  title   = {When and why {PINNs} fail to train: A neural tangent kernel perspective},
  journal = {Journal of Computational Physics},
  volume  = {449},
  pages   = {110768},
  year    = {2022},
  doi     = {10.1016/j.jcp.2021.110768}
}

@article{XPINNs,
  author  = {Jagtap, Ameya D. and Karniadakis, George Em},
  title   = {Extended physics-informed neural networks ({XPINNs}): A generalized space-time domain decomposition based deep learning framework for nonlinear partial differential equations},
  journal = {Communications in Computational Physics},
  volume  = {28},
  number  = {5},
  pages   = {2002--2041},
  year    = {2020},
  doi     = {10.4208/cicp.OA-2020-0164}
}

@article{cPINNs,
  author  = {Jagtap, Ameya D. and Kharazmi, Ehsan and Karniadakis, George Em},
  title   = {Conservative physics-informed neural networks on discrete domains for conservation laws: Applications to forward and inverse problems},
  journal = {Computer Methods in Applied Mechanics and Engineering},
  volume  = {365},
  pages   = {113028},
  year    = {2020},
  doi     = {10.1016/j.cma.2020.113028}
}

@article{Shin,
  author  = {Shin, Yeonjong and Darbon, Jerome and Karniadakis, George Em},
  title   = {On the convergence of physics informed neural networks for linear second-order elliptic and parabolic type {PDEs}},
  journal = {Communications in Computational Physics},
  volume  = {28},
  number  = {5},
  pages   = {2042--2074},
  year    = {2020},
  doi     = {10.4208/cicp.OA-2020-0193}
}

@article{RAR,
  author  = {Lu, Lu and Meng, Xuhui and Mao, Zhiping and Karniadakis, George Em},
  title   = {{DeepXDE}: A deep learning library for solving differential equations},
  journal = {SIAM Review},
  volume  = {63},
  number  = {1},
  pages   = {208--228},
  year    = {2021},
  doi     = {10.1137/19M1274067}
}

@article{RAD,
  author  = {Wu, Chenxi and Zhu, Min and Tan, Qinyang and Kartha, Yadhu and Lu, Lu},
  title   = {A comprehensive study of non-adaptive and residual-based adaptive sampling for physics-informed neural networks},
  journal = {Computer Methods in Applied Mechanics and Engineering},
  volume  = {403},
  pages   = {115671},
  year    = {2023},
  doi     = {10.1016/j.cma.2022.115671}
}

@article{FI_PINNs,
  author  = {Gao, Zhiwei and Yan, Liang and Zhou, Tao},
  title   = {Failure-informed adaptive sampling for {PINNs}},
  journal = {SIAM Journal on Scientific Computing},
  volume  = {45},
  number  = {4},
  pages   = {A1971--A1994},
  year    = {2023},
  doi     = {10.1137/22M1527763}
}

@article{DasPINNs,
  author  = {Tang, Kejun and Wan, Xiaoliang and Yang, Chao},
  title   = {{DAS-PINNs}: A deep adaptive sampling method for solving high-dimensional partial differential equations},
  journal = {Journal of Computational Physics},
  volume  = {476},
  pages   = {111868},
  year    = {2023},
  doi     = {10.1016/j.jcp.2022.111868}
}

@article{SA,
  author  = {McClenny, Levi D. and Braga-Neto, Ulisses M.},
  title   = {Self-adaptive physics-informed neural networks},
  journal = {Journal of Computational Physics},
  volume  = {474},
  pages   = {111722},
  year    = {2023},
  doi     = {10.1016/j.jcp.2022.111722}
}

\end{document}